\documentclass[oneside,english]{amsart}
\usepackage{bookman}
\usepackage[T1]{fontenc}
\usepackage[latin1]{inputenc}
\usepackage{amssymb}
\newcommand\Z{\ensuremath{\mathbb{Z}}}
\makeatletter

\providecommand{\LyX}{L\kern-.1667em\lower.25em\hbox{Y}\kern-.125emX\@}

 \theoremstyle{plain}
 \newtheorem*{lem*}{Lemma} 
 \theoremstyle{plain}
 \newtheorem*{prop*}{Proposition}
 \theoremstyle{remark}
 \newtheorem*{rem*}{Remark}
 \theoremstyle{plain}
 \newtheorem{thm}{Theorem}
 \newtheorem{prob}{Problem}

\def\th@nopoint{
  \thm@headpunct{} 
  \itshape 
}
\theoremstyle{nopoint}
\newtheorem*{ialt}{}

\usepackage{babel}
\makeatother
\begin{document}

\title{Random walks with $k$-wise independent increments}

\author{Itai Benjamini}

\author{Gady Kozma}

\author{Dan Romik}

\begin{abstract}
We construct examples of a random walk with pairwise-independent
steps which is almost-surely bounded, and for any $m$ and $k$ a
random walk with $k$-wise independent steps which has no
stationary distribution modulo  $m$.
\end{abstract}
\maketitle

\section{Introduction}

Consider a simulation of a simple random walk on a graph. How will
the simulation be affected if the source of randomness is not
truly random but only pseudo random in the following specific
sense, the random bits are $k$-wise independent for some $k>1$ and
not independent as a family? The first question to ask is, does the
pseudo walk converge to the same stationary measure? The most
simple graph to consider might be a cycle of size $m$. This
suggests the following problem: given a random walk $S_{n}=\sum
_{i=1}^{n}X_{i}$, where the $X_{i}$ are random signs, plus or
minus $1$ with equal probability, and the $X_{i}$'s are $k$-wise
independent, what can be said about the behavior of the partial
sums and in particular modulo some fixed number $m$? It turns out
that there is a fundamental difference between the cases $k=2$ and
$k>2$.

Examine first the case $k=2$. For this case we shall show

\begin{thm}
\label{thm:2bnd}There exists a sequence of random variables $\left\{ X_{i}\right\} _{i=1}^{\infty }$
taking values $\pm 1$, pairwise independent, such that $S_{n}$ is
bounded almost surely.
\end{thm}
This result should be contrasted against the fact that we know that
$\mathbb{E}S_{n}=0$ and $\mathbb{E}S_{n}^{2}=n$, just like in the
completely independent case. In other words, $S_{n}$ for large $n$
must have extreme {}``fat tail'' behavior. Naturally, $M:=\max _{n}S_{n}$
satisfies $\mathbb{E}M=\infty $. The example we will demonstrate
is essentially the Walsh system, which is pairwise independent. We
will discuss this in section \ref{sec:2in}.

Such behavior cannot occur when $k\geq 4$ for the simple reason
that in this case we know that $\mathbb{E}S_{n}^{4}=3n^{2}-2n$ and
this gives\[ \mathbb{P}(S_{n}^{2}>M)\geq
\frac{\left(\mathbb{E}S_{n}^{2}-M\right)^{2}}{\mathbb{E}S_{n}^{4}}\rightarrow
\frac{1}{3}\quad \forall M.\]

We could not settle $k =3$,
\begin{prob}
Is there   a sequence of random variables $\left\{ X_{i}\right\}
_{i=1}^{\infty }$ taking values $\pm 1$ with equal probability
$3$-wise independent, such that $S_{n}$ is bounded almost surely.
\end{prob}

The higher $k$ is, the more moments we know and we can approximate
the large scale shape of $S_{n}$.
However, this does not necessarily mean we can
say something about $S_{n}\mod m$. And indeed, in section
\ref{sec:Proofkm} we show the following

\begin{thm}
\label{thm:oneL}Let $m$ and $k$ be some natural numbers, and let
$\epsilon >0$. Then there exists a sequence of random variables $\left\{ X_{i}\right\} _{i=1}^{\infty }$
taking values $\pm 1$, $k$-wise independent, and a sequence $I_{j}$
such that \begin{equation}
\mathbb{P}\left(S_{I_{j}}\equiv 0\; (m)\right)>1-\epsilon \quad \forall j.\label{eq:thm}\end{equation}

\end{thm}
Notice that the requirement that the condition holds only for some
$I_{j}$ is unavoidable, since, say, if $k\geq 10m^{2}$ then the
distribution of $S_{I_{j}+10m^{2}}$ is approximately uniform, since
$X_{I_{j}+1},\dotsc ,X_{I_{j}+10m^{2}}$ are independent.
\medskip

Explicit constructions of $k$-wise independent 0-1 probability
spaces and estimates on their sizes are the focus of several
papers in combinatorics and complexity, see e.g.\ \cite{Jo} for the
first construction and \cite{AS} for a recent discussion with
additional references. Sums of pairwise independent random
variable were extensively studied, see e.g.\ \cite{Ja} for some
interesting examples. Thus there are already many interesting
examples of pairwise independent processes. But it seems behavior
modulo $m$ was not studied.

\section{\label{sec:2in}pairwise independent processes}

We term the construction we use a {}``gray walk'', as it is based
on the well-known Gray code construction in combinatorics. The $n$th
Gray code is a Hamiltonian path on the discrete $n$-cube, or a listing
of all $2^{n}$ binary strings of length $n$, starting with the string
$00...0$, where every string appears exactly once, and any two consecutive
strings differ in exactly one place. The construction (and hence also
proof that this is possible) is done recursively, namely: To construct
the $n$th Gray code, write down two copies of the $(n-1)$th code,
where the order of the strings in the second copy is reversed, and
add to the first $2^{n-1}$ strings a zero at the $n$th place, and
to the second $2^{n-1}$ strings a one at the $n$th place.

The Gray codes of all orders $n$ can be combined into an infinite
Gray code, by listing them sequentially for increasing $n$'s, and
converting all strings into infinite strings by padding with zeros.
The first few strings in the infinite code are:\begin{alignat*}{3}
A_{0} & =000000\ldots  & \qquad A_{1} & =100000\ldots  & \qquad A_{2} & =110000\ldots \\
A_{3} & =010000\ldots  & A_{4} & =011000\ldots  & A_{5} & =111000\ldots \\
A_{6} & =101000\ldots  & A_{7} & =001000\ldots  & A_{8} & =001100\ldots \\
 &  &  & \qquad \; \; \vdots  &  &
\end{alignat*}

To construct the gray walk, we now consider each string $A_{i}$ as
specifying a finite subset (also denoted $A_{i}$) of the natural
numbers $\mathbb{N}$, where a $1$ in the $j$th place signifies
that $j\in A_{i}$, and define \[
X_{i}=\prod _{j\in A_{i}}\xi _{j}\]
 where $\xi _{1},\xi _{2},...$ is a seed sequence of independent
$\pm 1$ variables. It is easy to verify that the $X_{i}$'s are pairwise
independent. Therefore to finish theorem \ref{thm:2bnd} we only need

\begin{prop*}
The gray walk $S_{n}$ is bounded almost surely. More precisely \[
\sup _{n}|S_{n}|=2^{\min \{j\geq 1:\xi _{j}=-1\}-1}+1=\sup S_{n}+2=-\inf _{n}S_{n}\]

\end{prop*}
\begin{proof}
For simplicity we add the element $X_{0}\equiv 1$, and prove that
for $S_{n}'=\sum _{i=0}^{n}X_{i}$ we have $\sup S_{n}'=-\inf S_{n}'=2^{\min \{j\geq 1:\xi _{j}=-1\}-1}$.
The recursive definition of the Gray code is reflected in the path:
Assume we have defined the first steps
$S_{0}',S_{1}',S_{2}',...,S_{2^{j-1}-1}'$
of the walk, then the next $2^{j-1}$ steps will be the previous steps
listed in reverse order, and multiplied by the random sign $\xi _{j}$.
This implies that the path up to time $2^{j}-1$, where $j$ is
the first value for which $\xi _{j}=-1$, is $0,1,2,3,4,...,2^{j-1},
2^{j-1}-1,2^{j-1}-2,...,3,2,1,0$,
and all the subsequent values are bounded between $-2^{j}$ and $2^{j}$.
\end{proof}
\begin{rem*}
Theorem \ref{thm:2bnd} holds also if one takes the strings $A_{i}$
in lexicographic order, which is simpler and is also the standard
ordering of the Walsh system. However, we believe the gray walk is
interesting in its own right. For example, it satisfies that
$\frac{X_{i+1}}{X_{i}}=\xi _{j(i)}$,
i.e.~the seed sequence is exposed directly in the quotients of the
$X_{i}$'s.

On the other hand, there is a vast body of knowledge about the behavior of
``weighted walks'' with the lexicographical order, which is simply the question
when
\[
\sum_{i=0}^{\infty} c_{i} W_{i}
\]
converges where $W_{i}$ is the Walsh system. The general philosophy is that it
is similar to the behavior of the Fourier system, with the exception of the
symmetry of the Fourier system which has no equivalent in the Walsh system. Of
course, this is {\it very different} from the simple behavior of independent
variables. See \cite{W97} for a survey, or the book \cite{SWSP90}.
\end{rem*}

\section{\label{sec:Proofkm}Proof of theorem \ref{thm:oneL}}

Before embarking on the general proof, let us demonstrate the case
$m=4$ which is far simpler. Let $L>k$ satisfy $L\equiv 0$ mod $4$.
Let $\left\{ \xi _{i}\right\} _{i=1}^{\infty }$ be a sequence of
i.i.d $\pm 1$ variables. Next define new variables by conditioning the $\xi_i$-s:\[
\left\{ X_{i}\right\} :=\left\{ \xi _{i}\left|\prod _{b=1}^{L}\xi _{aL+b}=1\; \forall a=0,1,\dotsc \right.\right\} \quad .\]
The fact that $L>k$ clearly shows that the $X_{i}$'s are $k$-independent.
The fact that $\prod \xi _{aL+b}=1$ shows that in each block of size
$L$ the number of $-1$'s is even, and since $4|L$ we get that $\sum \xi _{aL+b}\equiv 0$
mod 4 for every block. Therefore if we define $I_{j}=jL$ then the
condition in (\ref{eq:thm}) is actually satisfied combinatorially
and not only in high probability.

How can we generalize this to $m\neq 4$? The remaining ingredient
in the proof is based on the following well known fact:

\begin{ialt}It is possible to construct variables with probability
$\frac{1}{3}$, $\frac{1}{3}$, $\frac{1}{3}$ given a sequence of
independent variables with probabilities $\frac{1}{2}$, $\frac{1}{2}$.\end{ialt}

The algorithm is simple: take two such variables: if they give $11,$
{}``output'' $1$, if they give $1,-1$ {}``output'' 2 and if
they give $-1,1$, {}``output'' $3$. If they give $-1,-1$, take
a new set of two variables and repeat the process. The output is $1$,
$2$ or $3$, each with probability $\frac{1}{3}$. The time required
for each step is unbounded but the probability that it is large is
exponentially small.

The proof below combines these two ideas, which we nickname {}``generating
hidden symmetries'' and {}``simulating uniform variables'' to get
the result.

\begin{proof}[Proof of theorem \ref{thm:oneL}, general $m$]We may assume without
loss of generality that $m$ is even (otherwise take $2m$ instead).
Let $\lambda =\lambda (k,\epsilon )$ be some sufficiently large parameter
which we shall fix later. Let $L:=2(k+1)\left\lceil \lambda m^{2}\right\rceil $
where $\lceil \cdot \rceil$ stands, as usual, for the {\it upper} integer
value.

\begin{lem*}
For every $\mu =0,2,\dotsc ,m-2$ there exists variables $L^{\mu }$
and $Y^{\mu }$ with the following properties
\begin{enumerate}
\item \label{enu:LmuL}$L^{\mu }$ is a random positive integer, $L$ divides
$L^{\mu }$ and $\mathbb{P}(L^{\mu }=aL)\leq 2^{-a-1}\epsilon $ for
all $a>1$.
\item \label{enu:Xmukindep}$Y^{\mu }=\left\{ y_{i}^{\mu }\right\} _{i=1}^{\infty }$
is a sequence of $k$-wise independent variables taking values $\pm 1$
with probabilities $\frac{1}{2}$.
\item \label{enu:Xmularge}$\left\{ \left.\left\{ y_{i\textrm{ }}^{\mu }\right\} _{i=L^{\mu }+1}^{\infty }\right|L^{\mu }\right\} $
are i.i.d $\pm 1$ variables
\item \label{enu:sumxmu}$L^{\mu }=L$ implies $\sum _{i=1}^{L}y_{i}^{\mu }\equiv \mu $
mod $m$.
\end{enumerate}
\end{lem*}
\begin{proof}
Define $N:=2\left\lceil \lambda m^{2}\right\rceil $. The first step is to
divide $\{\pm 1\}^N$ according to the sum modulo $m$, and trim the resulting
sets a little so that they all have the same size. Precisely, let\[
A:=\min _{j =0,2,\dotsc ,m-2}\#\left\{ v\in \{\pm 1\}^{N}:\sum _{i=1}^{N}v_{i}\equiv j \; (m)\right\} \quad .\]
We note that the distribution of $\sum _{i=1}^{N}\pm 1$ modulo $m$
is uniform on the set $0,2,\dotsc ,m-2$ with an error of $Ce^{-c\lambda }$.
Therefore $2^{N}-\frac{m}{2}A\leq C2^{N-c\lambda }$. For
$j =0,2,\dotsc ,m-2$,
let $G_{j}$ be arbitrary sets satisfying \[
G_{j}\subset \left\{ v\in \{\pm 1\}^{N}:\sum _{i=1}^{N}v_{i}\equiv \lambda \; (m)\right\} \quad |G_{j}|=A.\]
and let $S:=\{\pm 1\}^{N}\setminus \bigcup G_{j}$.

So far we have constructed some general objects having very little to
do with probability. Next, let $\left\{ \xi _{i}\right\} _{i=1}^{\infty }$
be i.i.d variables taking values in $\pm 1$, and let $\left\{ \Xi _{i}\right\} _{i=1}^{\infty }$
be a division of the $\xi _{i}$'s into blocks of size $N$: $\Xi _{i}:=(\xi _{Ni-N+1},\dotsc ,\xi _{Ni})$.
Let \[
B:=\min \left\{ b:\#\left\{ \Xi _{i}\not \in S\right\} _{i=1}^{b}=k+1\right\} \quad .\]
We define\begin{align*}
L^{\mu } & :=\left\{ L\left\lceil \frac{B}{k+1}\right\rceil \left|\sum _{i=1}^{BN}\xi _{i}\equiv \mu \; (m)\right.\right\} \\
Y^{\mu } & :=\left\{ (\xi _{1},\xi _{2},\dotsc )\left|\sum _{i=1}^{BN}\xi _{i}\equiv \mu \; (m)\right.\right\} \quad .
\end{align*}

Properties \ref{enu:Xmularge} and \ref{enu:sumxmu} are obvious from
the construction. Property \ref{enu:LmuL} is a direct consequence
of the estimate $\mathbb{P}(\Xi _{i}\in S)\leq Ce^{-c\lambda }$,
if only $\lambda $ is large enough. Define $\lambda $ so as to satisfy
this condition. Therefore we need only prove property \ref{enu:Xmukindep}.

Let therefore $i_{1}<\dotsb <i_{k}$ and $\delta _{1},\dotsc ,\delta _{k}\in \{\pm 1\}$
and examine the events \[
\mathcal{X}:=\left\{ \xi _{i_{1}}=\delta _{1},\dotsc ,\xi _{i_{k}}=\delta _{k}\right\} \quad \mathcal{Y}:=\left\{ \sum _{i=1}^{BN}\xi _{i}\equiv \mu \; (m)\right\} .\]
We know that $\mathbb{P}(\mathcal{X})=2^{-k}$ and we need to show
that $\mathbb{P}(\mathcal{X}|\mathcal{Y})=2^{-k}$. Let $b$ be some
number and let $\sigma \subset \{1,\dotsc ,b-1\}$ be a set with $\#\sigma =b-(k+1)$
and examine the event \[
\mathcal{S}=\mathcal{S}(\sigma ,b)=\{\xi _{i}\in S\Leftrightarrow i\in \sigma \quad \forall i\leq b\}\quad .\]
The point of the proof is that the event $\mathcal{Y}$ is independent
of $\mathcal{X}\cap\mathcal{S}$. This is because $\mathcal{X}$
depends only on $k$ places, but $\mathcal{Y}$ depends on $k+1$
$\Xi _{i}$'s each of which is distributed uniformly. In other words\[
\mathbb{P}(\mathcal{Y}|\mathcal{X}\cap \mathcal{S})=\mathbb{P}(\mathcal{Y})=\frac{2}{m}\]
or \[
\mathbb{P}(\mathcal{X}\cap \mathcal{S}|\mathcal{Y})=\frac{\mathbb{P}(\mathcal{X}\cap \mathcal{S}\cap \mathcal{Y})}{\mathbb{P}(\mathcal{Y})}=\mathbb{P}(\mathcal{X}\cap \mathcal{S})\]
which gives \[
\mathbb{P}(\mathcal{X}|\mathcal{Y})=\sum _{\sigma ,b}\mathbb{P}(\mathcal{X}\cap S|\mathcal{Y})=\sum _{\sigma ,b}\mathbb{P}(\mathcal{X}\cap \mathcal{S})=\mathbb{P}(\mathcal{X})=2^{-k}\quad .\qedhere \]

\end{proof}
Returning to the proof of the theorem, we let \[
\left\{ Y^{\mu ,\nu },L^{\mu ,\nu }\right\} _{(\mu ,\nu )\in \{0,2,\dotsc ,m-2\}\times \mathbb{N}}\]
 be an independent family of couples of variables constructed using
the lemma. We now construct our sequence $X_{i}$ inductively as follows:
assume that at the $\nu $th step we have defined $X_{1},\dotsc ,X_{\rho}$
where $\rho=\rho(\nu )$ is random. Define
$\mu =\mu (\nu )\equiv -\sum _{i=1}^{\rho}X_{i}$
mod $m$ and then define $X_{\rho+1},\dotsc ,X_{\rho+L^{\mu ,\nu }}$ using
\[
X_{\rho+i}=y_{i}^{\mu ,\nu }\]
and $\rho(\nu +1)=\rho(\nu)+L^{\mu ,\nu }$. This creates a sequence of $\pm 1$
variables. To see (\ref{eq:thm}), define $r(a)=\max \{\nu :\rho(\nu )<aL\}$
and properties \ref{enu:LmuL} and \ref{enu:sumxmu} of the lemma
will give that
\begin{align*}
\mathbb{P}\left(\sum _{i=1}^{aL}\epsilon _{i}\not \equiv 0\; (m)\right) &
\leq \sum _{b=1}^{\infty }\mathbb{P}\big(\rho(r(a))=(a-b)L\textrm{ and }
L^{\mu (r(a)),r(a)}\geq L\max\{b,2\}\big)\\
 & \leq \sum _{b=1}^{\infty }\sum _{c=\max\{b,2\}}^{\infty }
\mathbb{P}(L^{\mu (r(a)),r(a)}=cL)\leq \sum _{b,c}2^{-c-1}\epsilon =
\epsilon \quad .
\end{align*}
Therefore we need only show that the $X_{i}$'s are $k$-wise independent.
While being strongly related to the $k$-wise independence of the
$y_{i}^{\mu ,\nu }$
it is not an immediate consequence of it because of the dependence
between the $y_{i}^{\mu ,\nu }$ and the $L^{\mu ,\nu }$'s.

We prove that the $X_{i}$'s are $k$-wise independent inductively. Let
$l$ and $I$ be some integers and assume that we have shown that
$X_{i_{1}},\dotsc ,X_{i_{m}}$ are independent for every $m<l$ and
for $m=l$ and $i_{1}<I$. Let $I=i_{1}<i_{2}<\dotsb <i_{l}$ and
let $\delta _{1},\dotsc ,\delta _{l}\in \{\pm 1\}$. Examine $L^{0,1}$
and define\[
n=n(L^{0,1})=\#\{m:i_{m}\leq L^{0,1}\}\quad .\]
The induction hypothesis, together with the independence of the various
$Y^{\mu ,\nu }$'s shows that $X_{i_{n+1}},\dotsc ,X_{i_{l}}$ are
i.i.d $\pm 1$ variables: if $n>0$ then this follows from the induction
hypothesis for $m':=m-n$ while if $n=0$ then it follows from the induction
hypothesis for $m'=m$ and $I':=I-L^{0,1}$. Property \ref{enu:Xmularge} of the lemma
shows that $y_{i_{n+1}}^{0,1},\dotsc ,y_{i_{l}}^{0,1}$ are i.i.d
$\pm 1$ variables. Below $n$ we have simple equality:\[
X_{i_{m}}=y_{i_{m}}^{0,1}\quad \forall m\leq n.\]
Therefore\[
\mathbb{P}(X_{i_{1}}=\delta _{1},\dotsc ,X_{i_{l}}=\delta _{l}|L^{0,1})=\mathbb{P}(y_{i_{1}}=\delta _{1},\dotsc ,y_{i_{l}}=\delta _{l}|L^{0,1})\]
 and summing over $L^{0,1}$ we get the required result:\[
\mathbb{P}(X_{i_{1}}=\delta _{1},\dotsc ,X_{i_{l}}=\delta _{l})=\mathbb{P}(y_{i_{1}}=\delta _{1},\dotsc ,y_{i_{l}}=\delta _{l})=2^{-l}\]
which finishes the proof of the theorem.\end{proof}

\theoremstyle{remark}

\newtheorem*{rmks}{Remarks}

\begin{rmks}

\begin{enumerate}
\item It is easy to generalize the theorem to force $\sum _{i=1}^{I_{j}}\epsilon _{i}$
to have almost arbitrary distributions. The only obstacle is the even-odd
symmetry of the walk. Thus if $m$ is odd, any distribution whatsoever
on the congruence classes modulus $m$ might be achieved, while if
$m$ is even, any distribution supported on the set of even or odd
congruence classes can be achieved.
\item It is easy to see from the proof that $L\approx km^{2}\log (k/\epsilon )$.
In other words, the precision is asymptotically exponential in the
stepping of the subsequence $I_{j}$. We remind the reader a fact
that was discussed in the introduction: if $k>m^{2}\log 1/\epsilon $
then necessarily $I_{j+1}-I_{j}>cm^{2}\log 1/\epsilon $ because immediately
after $I_{j}$ we get a sequence of independent variables.

\item

Another direction in which $k$-wise independent random variables
may be very different from truly independent variables is their
entropy. It is known that the entropy of $n$ such variables may be
as low as $k\log n$, see \cite{Jo} and \cite{AS}. Our pairwise
example is optimal in this respect: $n$ truly random bits create
$2^n$ pairwise independent random variables. The example above is
far from it: the entropy is linear in the number of bits. However,
it turns out that it is possible to create an example of $k$-wise
independent variables with low entropy satisfying the conditions
of theorem \ref{thm:oneL}. The example, while amusing, if a bit
off topic hence we will skip it.

\end{enumerate}
\end{rmks}

\begin{thm}
Let $m$ and $k$ be some natural numbers. Then there exists a sequence
of random variables $\left\{ \epsilon _{i}\right\} _{i=1}^{\infty }$
taking values $\pm 1$, $k$-wise independent, and a sequence $I_{j}$
such that \[
\sum _{i=1}^{I_{j}}\epsilon _{i}\equiv 0\; (m)\textrm{ for }j\textrm{ sufficiently large with probability }1.\]

\end{thm}
The proof is similar to that of theorem \ref{thm:oneL}, but using
a different $L$ in each step. However there are some additional technical
subtleties. Essentially, if in the main lemma of the previous theorem
we defined a variable $Y^{\mu }$ where the parameter $\mu $ was
used to {}``return to sync'' the series if ever the congruence is
no longer $0$, here we would need variables $Y^{\mu ,\tau }$ where
$\tau $ is some additional parameter needed to {}``return to sync''
in the $L$ domain. While being only moderately more complicated,
we felt it better to present the simpler proof of the previous theorem.

\section{Further related problems}
We end by suggesting further problems related to $k$-wise
independent random variables.

\subsection{The random sign Pascal triangle}

Let $\xi _{n,k}^{+},\xi _{n,k}^{-}$ be a family of i.i.d random
(1/2-1/2) signs for $-1\leq k\leq n+1$. Define random variables
$X_{n,k}$ for $-1\leq k\leq n+1$ by the recurrence \[
X_{0,0}=1,\quad X_{n,-1}=0,\quad X_{n,n+1}=0, \]
 \[
X_{{n,k}}=\xi_{n-1,k-1}^{+}X_{n-1,k-1}+\xi_{n-1,k}^{-}X_{n-1,k}
\quad (0\le k\le n)
\]

Problems: Study the behavior of the central coefficient
$X_{2n,n}$. Find an interpretation of the triangle as {}``percolation
with interference'' where two closed gates cancel each other out
({}``quantum percolation''?) Study the behavior of the non-central
coefficients.

\subsection{2D Percolation}

What about other functions of $k$-wise independent random
variables? For example, instead of considering Bernoulli
percolation, the random variables determining the configuration
are only $k$-wise independent. Here is an  ``unimpressive''
example:
  $k$-wise independent bond percolation on an $n  \times n(k+1)$- box
  in $\Z^2$ such
that the probability to have a crossing of the narrow side is $1$.
Of course, even in usual percolation the
probability\footnote{This is obvious if you are willing to use the conformal
invariance of percolation. Otherwise, this fact can be derived from
Russo-Seymour-Welsh theory.} is $1-e^{-ck}$
so that the improvement is not exciting.

The construction is as follows: divide into $k+1$ disjoint boxes of size
$n  \times n$. Take the usual independent percolation, conditioned such
that the number of boxes with crossings is odd. This is like
conditioning $k+1$ Bernoulli variables to have an odd sum, so that
you get a $k$-wise independent structure. If the number is odd,
then it is non zero. Can it be improved, say to guarantee crossing
of the $n \times n$-box?

\end{document}